\documentclass[twoside,12pt]{article}
\usepackage{amssymb,amsmath,bm,euscript}
\usepackage{graphicx,color,caption}
\usepackage{ifpdf}

\usepackage[colorlinks=true]{hyperref}

\usepackage{braket}

\usepackage{algorithm, algorithmic}

\newtheorem{proposition}{Proposition}[section]
\newtheorem{theorem}[proposition]{Theorem}

\newtheorem{corollary}[proposition]{Corollary}
\newtheorem{definition}[proposition]{Definition}

\newcommand{\qed}{\hphantom{.}\hfill $\Box$\medbreak}

\def\I{\mathcal{I}}

\def\A{{\mathcal{A}}}

\def\D{\mathcal{D}}

\def\X{{\mathcal{X}}}
\def\Y{{\mathcal{Y}}}
\def\Z{\mathcal{Z}}

\def\M{\mathcal{M}}

\def\aa{{\bf a}}
\def\bb{{\bf b}}

\def\ee{{\bf e}}
\def\x{{\bf x}}

\def\0{{\bf 0}}

\title{\bf{A Tensor Rank Theory and Maximum Full Rank Subtensors}}%\thanks{This research was supported by the Hong Kong
\author{ \hspace{1mm} Liqun Qi\thanks{%Department of Mathematics, School of Science, Hangzhou Dianzi University, Hangzhou 310018 China; Future Network Theory Lab, 2012 Labs
Department of Applied
    Mathematics, The Hong Kong Polytechnic University, Hung Hom,
    Kowloon, Hong Kong, China; ({\tt liqun.qi@polyu.edu.hk}).},
 \ \
  Xinzhen Zhang\thanks{School of Mathematics, Tianjin University, Tianjin 300354 China; ({\tt xzzhang@tju.edu.cn}). This author's work was supported by NSFC (Grant No.  11871369). },
 \ and \
Yannan Chen\thanks{School of Mathematical Sciences, South China    Normal University, Guangzhou, China; ({\tt ynchen@scnu.edu.cn}).  This author was supported by the National Natural Science Foundation of China (11771405).}
}

\begin{document}
\date{\today}
\maketitle

\begin{abstract}
A matrix always has a full rank submatrix such that the rank of this matrix is equal to the rank of that submatrix.  This property is one of the corner stones of the matrix rank theory.  We call this property the max-full-rank-submatrix property.    Tensor ranks play a crucial role in low rank tensor approximation, tensor completion and tensor recovery.   However, their theory is still not matured yet.   Can we set an axiom system for tensor ranks?   Can we extend the max-full-rank-submatrix property to tensors?   We explore these in this paper. We first propose some axioms for tensor rank functions.   Then we introduce proper tensor rank functions.  The CP rank is a tensor rank function, but is not proper.   There are two proper tensor rank functions, the max-Tucker rank and the submax-Tucker rank, which are associated with the Tucker decomposition.  We define a partial order among tensor rank functions and show that there exists a unique smallest tensor rank function.   We introduce the full rank tensor concept, and define the max-full-rank-subtensor property.   We show the max-Tucker tensor rank function and the smallest tensor rank function have this property.    We define the closure for an
arbitrary proper tensor rank function, and show that it is still a proper tensor rank function and has the max-full-rank-subtensor property.   An application of the submax-Tucker rank is also presented.

%Our theoretic analysis indicates that the submax-Tucker rank is a good choice for low rank tensor approximation and tensor completion.


%15A69:LINEAR AND MULTILINEAR ALGEBRA; MATRIX THEORY- Multilinear algebra, tensor products
%53A45:Classical differential geometry-Vector and tensor analysis
%47A05: Operator Theory-General (adjoints, conjugates, products, inverses, domains, ranges,
%etc.)
%53C35:Global differential geometry- Symmetric spaces
\end{abstract}

%\newpage

\section{Introduction}

A matrix always has a full rank submatrix such that the rank of this matrix is equal to the rank of that submatrix.   We call this property the max-full-rank-submatrix property.
This property is one of the corner stones of the matrix rank theory.

We now arrive the era of big data and tensors.  Tensor ranks play a crucial role in low rank tensor approximation, tensor completion and tensor recovery \cite{ADKM11, JYZ17, KBHH13, KB09, TYFWR13, XWWXWZ16, YHHH16, ZSKA18, ZA17, ZEAHK14, ZZXC15, ZLLZ18}.   However, their theory is still not matured yet.   Can we
set an axiom system for tensor ranks?   Can we extend the full rank concept and the max-full-rank-submatrix property to tensors?   We explore these in this paper.

We first propose some axioms for tensor rank functions.   Then we introduce proper tensor rank functions.  The CP rank is a tensor rank function, but is not proper.   There are two proper tensor rank functions, the max-Tucker rank and the submax-Tucker rank, which are associated with the Tucker decomposition.  We define a partial order among tensor rank functions and show that there exists a unique smallest tensor rank function.   We introduce the full rank tensor concept, and define the max-full-rank-subtensor property.   We show the max-Tucker tensor rank function and the smallest tensor rank function have this property.    We define the closure for an
arbitrary proper tensor rank function, and show that it is still a proper tensor rank function and has the max-full-rank-subtensor property.   An application of the submax-Tucker rank is also presented.

The set of all nonnegative integers is denoted by ${\bf Z_+}$.   The set of all positive integers is denoted by ${\bf N}$.    Let $m, n_1, \cdots, n_m \in {\bf N}$.
Denote the set of all real $m$th order tensors of dimension $n_1 \times n_2 \times \cdots \times n_m$ by ${\bf T}(n_1, n_2, \cdots, n_m)$.     If $n_1 = \cdots = n_m = n$, then we denote it by ${\bf CT}(m, n)$.  Here ``CT'' means cubic tensors.   Denote the set of all real tensors by ${\bf T}$.   Thus, scalars, vectors, matrices are a part of  ${\bf T}$.   Let $\X \in {\bf T}(n_1, n_2, \cdots, n_m)$.   We call
$\X$ a rank-one tensor if and only if there are nonzero vectors $\x^{(i)} \in \Re^{n_i}$ for $i = 1, \cdots, m$, such that
$$\X = \x^{(1)} \circ \cdots \circ \x^{(m)}.$$
Here, $\circ$ is the tensor outer product.   Then, nonzero vectors and scalars are all rank-one tensors in this sense.

Suppose that $m, n_1, \cdots, n_m \in {\bf N}$ and $\X = (x_{i_1\cdots i_m}) \in {\bf T}(n_1, \cdots, n_m)$.  Let $T_l \subset \{ 1, \cdots, n_l \}$ and $|T_l | = k_l \ge 1$ for $l = 1, \cdots, m$.   Suppose that $\Y = (y_{j_1\cdots j_m}) \in {\bf T}(k_1, \cdots, k_m)$ for $j_l \in T_l$, with $y_{j_1\cdots j_m} = x_{i_1\cdots i_m}$ if $j_l=i_l$ for $l = 1, \cdots, m$.  Then we say that $\Y$ is a subtensor of $\X$.    If $\Y \not = \X$, then we say that $\Y$ is a proper subtensor of $\X$.  %If all the entries of $\X$, which are not in $\Y$, are zero, then we say that $\Y$ is an essential subtensor of $\X$.
For $p = 1, \cdots, m$, the subtensor $\Y$, described above, is called a $p$-row of $\X$ if $|T_p| = 1$ and $T_l = \{ 1, \cdots, n_l \}$ for $l \not = p$.    If $T_p = \{ q \}$,  then the corresponding $p$-row is called the $q$th $p$-row of $\X$.    The $p$-row concept extends the concepts of rows and columns from matrices to tensors.  For a matrix, a $1$-row is called a row, a $2$-row is called a column.

In the next section, we present a set of axioms for tensor rank functions.   We list six properties which are essential for tensor ranks.   In particular, we define a partial order ``$\le$'' among tensor rank functions,
and show that there exists a unique smallest tensor rank function $r_*$.  We also introduce proper and strongly proper tensor rank functions in that section.

 %It is interesting to identify this smallest tensor rank function $r_*$ and its properties.    %This is significant in theory and may have potential application values.

We study the CP rank and the Tucker rank in Section 3.    The Tucker rank is a vector rank.   We derive two scalar ranks from this, and call them the max-Tucker rank and the submax-Tucker rank respectively.   We show that the CP rank, the max-Tucker rank and the submax-Tucker rank are all tensor rank functions.   %The max-Tucker rank naturally arises from the Tucker decomposition.   The submax-Tucker rank is also associated with the Tucker decomposition, but is a new tensor rank function introduced in this paper.
The CP rank is subadditive but not proper.   The max-Tucker rank is proper, subadditive but not strongly proper.  The submax-Tucker rank is strongly proper but not subadditive.

%The rank of a matrix is always equal to the maximum value of the ranks of all its full rank submatrices.   We call this property the sub-full-rank property.   It forms a cornerstone of the matrix rank theory.

We introduce the concept of maximum full rank subtensors, and define the max-full-rank-subtensor property
in Section 4.   We show that the max-Tucker rank function has the max-full-rank-subtensor property.     Suppose that $m, n_1, \cdots, n_m \in {\bf N}$ and $\X = (x_{i_1\cdots i_m}) \in {\bf T}(n_1, \cdots, n_m)$, and $\Y$ is a maximum full rank subtensor of $\X$ under the max-Tucker rank.
Then we show that there is an index $p$, $1 \le p \le m$, such that all the $q$th $p$-rows of $\X$, with $q$ in the mode $p$ index set $T_p$ of $\Y$, are linearly independent, and any $p$-row of $\X$ is a linear combination of the $q$th $p$-rows of $\X$ with $q \in T_p$.

In Section 5, we define the closure of an arbitrary proper tensor rank function, and show that it is still a proper tensor rank function, and has the max-full-rank-subtensor property.   We show that $r_*$ is strongly proper and has the max-full-rank-subtensor property.

We present an application of the submax-Tucker rank in internet traffic data approximation in Section 6.

Some final remarks are made in Section 7.

We use small letters to denote scalars, small bold letters to denote vectors, capital letters to denote matrices, and calligraphic letters to denote tensors.   We denote the matrix rank of a matrix $A$ as $r_0(A)$.

\section{Axioms and Properties of Tensor Rank Functions}

Let $m, n \in {\bf N}$.   Consider ${\bf CT}(m, n)$.   Suppose $\X = (x_{i_1\cdots i_m}) \in {\bf CT}(m, n)$.  An entry $x_{i_1\cdots i_m}$ is called a diagonal entry of $\X$ if $i=i_1=\cdots = i_m$.   Otherwise, $x_{i_1\cdots i_m}$ is called an off-diagonal entry of $\X$.   If all the off-diagonal entries of $\X$ are zero, then $\X$ is called a diagonal tensor.  If $\X \in {\bf CT}(m, n)$ is diagonal, and all the diagonal entries of $\X$ is $1$, then $\X$ is called the identity tensor of ${\bf CT}(m, n)$, and denoted as $\I_{m, n}$.   Clearly, the identity tensor $\I_{m, n}$ is unique to ${\bf CT}(m, n)$.   The identity tensor plays an important role in spectral theory of tensors \cite{QL17}.

\begin{definition} \label{d0}
Suppose that $r : {\bf T} \to {\bf Z_+}$.   If $r$ satisfies the following six properties, then $r$ is called a tensor rank function.

{\bf Property 1}  Suppose that $\X \in {\bf T}$.   Then $r(\X) = 0$ if and only if $\X$ is a zero tensor, and
$r(\X) =1$ if and only if $\X$ is a rank-one tensor.

{\bf Property 2} For $m, n \in {\bf N}$ with $m \ge 2$, $r(\I_{m, n}) = n$.

{\bf Property 3}  Let $n_1, n_2 \in {\bf N}$, $\X \in {\bf T}(n_1, n_2, 1, 1, \cdots, 1)$.   Then $r(\X)$ is equal to the matrix rank of the $n_1 \times n_2$ matrix corresponding to $\X$.

{\bf Property 4}  Let $m, n_1, \cdots, n_m \in {\bf N}$, $\X = (x_{i_1\cdots i_m}) \in {\bf T}(n_1, \cdots, n_m)$, and $\alpha$ is a real nonzero number.    Then $r(\X) = r(\alpha\X)$.

{\bf Property 5}  Let $m, n_1, \cdots, n_m \in {\bf N}$, $\X = (x_{i_1\cdots i_m}) \in {\bf T}(n_1, \cdots, n_m)$, and $\sigma$ is a permutation on ${\bf N}^m$.    Then $r(\X) = r(\Y)$, where $\Y = (x_{j_1\cdots j_m}) \in
{\bf T}(\sigma(n_1, \cdots, n_m))$, $(j_1, \cdots, j_m) = \sigma(i_1, \cdots, i_m)$.

{\bf Property 6}   Let $m, n_1, \cdots, n_m \in {\bf N}$.   Suppose that $\X, \Y \in {\bf T}(n_1, \cdots, n_m)$, and $\Y$ is a subtensor of $\X$.  Then $r(\Y) \le r(\X)$.   %If $\Y$ is an essential subtensor of $\X$, then $r(\Y) = r(\X)$.
\end{definition}

These six properties are essential for tensor ranks.  Property 1 specifies rank zero tensors and rank-one tensors.   Though the tensor rank theory is not matured, there are no arguments in rank zero and rank-one tensors in the literature.   Property 2 fixes the value of the tensor rank for identity tensors.  This is necessary as identity tensors are good references for the magnitude of tensor ranks.   Property 3 justifies the tensor rank is an extension of the matrix rank.   Property 4 claims that the tensor rank is not changed when a tensor is multiplied by a nonzero real number.   Property 5 says that the roles of the modes are balanced.   Property 6 justifies the subtensor rank relation.

Suppose that $r_1, r_2 : {\bf T} \to {\bf Z_+}$ are two tensor rank functions.
If for any $\X \in {\bf T}$ we always have $r_1(\X) \le r_2(\X)$, then we say that the tensor rank function
$r_1$ is not greater than the tensor rank function $r_2$ and denote this relation as $r_1 \le r_2$.

\begin{theorem}  \label{minth}
Suppose that $r_1, r_2 : {\bf T} \to {\bf Z_+}$ are two tensor rank functions.  Define
$r: {\bf T} \to {\bf Z_+}$ by
$$r(\X) = \min \{ r_1(\X), r_2(\X) \},$$
for any $\X \in {\bf T}$.  Then $r$ is a tensor rank function, $r \le r_1$ and $r \le r_2$.
\end{theorem}
{\bf Proof}   For any $\X\in T$, let $r(\X)=\min \{ r_1(\X), r_2(\X) \}$.  Then Properties 1, 2, 3 and 4 hold clearly from the definition of tensor rank functions.

To show Property 5, we assume that $\Y$ is a permutated tensor of $\X$. Then $r_1(\X)=r_1(\Y)$ and $r_2(\X)=r_2(\Y)$. Hence, $r(\X)=r(\Y)$ and Property 5  is obtained.

Now we assume that $\Z$ is a subtensor of $\X\in {\bf T}$. Then $r_1(\Z)\leq r_1(\X)$ and $r_2(\Z)\leq r_2(\X)$.
Hence $r(\Z)\leq r(\X)$ since $r(\Z)=\min \{r_1(\Z),r_2(\Z)\}\leq r_1(\X)$ and $r(\Z)\leq r_2(\X)$. %For an essential subtensor $\Z$, $r_1(\X)=r_1(\Z)$ and $r_2(\X)=r_2(\Z)$. Hence $\min\{r_1(\Z),r_2(\Z)\}=\min\{r_1(\X), r_2(\X)\}$ and
Thus, Property 6 holds.

Thus, we conclude that $r=\min \{r_1, r_2\}$ is a tensor rank function.

clearly, $r \le r_1$ and $r \le r_2$.

\qed

\begin{theorem}
There exists a unique tensor rank function $r_*$, such that for any tensor rank function $r$, we have $r_* \le r$.
\end{theorem}
{\bf Proof}  For any $\X \in {\bf T}$, define $r_*(\X):= \min \{r(\X) \ | \ r \mbox{ is~~ a ~~tensor~~ rank ~~ function}\}$.    This is well-defined as tensor rank functions take values on $\bf Z_+$.  Now we show that $r^*$ is a tensor rank function.

1) Suppose $\X$ is a zero tensor in $\bf T$.   Then for any tensor rank function $r$, $r(\X) = 0$.  This implies that $r_*(\X) = 0$ by the definition of $r_*$.   On the other hand, suppose that $r_*(\X)=0$ for some $\X \in {\bf T}$.  Then for some tensor rank function $r$, $r_*(\X)=r(\X)=0$. Hence, $\X$ is a zero tensor from Property 1 of the tensor rank function $r$.   Similarly, we may show that $r_*(\X) = 1$ if and only if $\X$ is a rank-one tensor.

2) For any $m, n \in {\bf N}$ with $m \ge 2$, $r(\I_{m,n})=n$ for all tensor rank functions $r$. Thus $r_*(\I_{m, n})=n$.

 3) Let $\X\in {\bf T}(n_1,n_2,1,\cdots, 1)$. Let $M$ be the corresponding $n_1 \times n_2$ matrix in $\X$.  Then for any tensor rank function $r$, $r(\X) = r_0(M)$.   Hence all of $r(\X)$ are equal. Hence, $r_*(\X) = r_0(M)$ and Property 3 holds.

4) For any $\X\in {\bf T}$ and any tensor rank function $r$, $r(\X)=r(\alpha \X)$ for any $\alpha\neq 0$. Thus, $r_*(\X)=r_*(\alpha\X)$.

5) We have Properties 5 and 6 in a similar way as in the proof of Theorem \ref{minth} and omit the details here.

By the definition, $r_* \le r$ for any tensor rank function $r$.

Suppose that $r_*$ and $r_{**}$ are two tensor rank functions with the property that $r_* \le r$ and $r_{**} \le r$ for any tensor rank function.   Then $r_* \le r_{**} \le r_*$.  We see that $r_* = r_{**}$.   Thus, such a tensor rank function $r_*$ is unique.

\qed

We call $r_*$ the smallest tensor rank function.   In Section 5, we will show that $r_*$ has the max-full-rank-subtensor property.

The six properties in Definition \ref{d0} are essential to tensor rank functions.
There are some other properties which are satisfied by some tensor rank functions.

\begin{definition}
Suppose that $r$ is a tensor rank function.  We say that $r$ is a proper tensor rank function if for any $m, n \in {\bf N}$ and $\X \in {\bf CT}(m, n)$, we have $r(\X) \le n$.
\end{definition}

For an $n \times n$ square matrix, its matrix rank is never greater than its dimension $n$.
Thus, proper tensor rank functions are reasonable in a certain sense.

\begin{definition}
Suppose that $r$ is a tensor rank function.  We say that $r$ is a subadditive tensor rank function if for any $m, n_1, \cdots, n_m \in {\bf N}$, and $\X, \Y \in {\bf T}(n_1, \cdots, n_m)$, we have
$$r(\X + \Y) \le r(\X) + r(\Y).$$
\end{definition}

The subadditivity property is somehow restrictive.   The minimum of two subadditive tensor rank functions may not be subadditive.

\begin{proposition}
Suppose that $r$ is a proper tensor rank function.   Let $m, n_1, \cdots, n_m \in {\bf N}$ with $m \ge 2$ and $\X \in {\bf T}(n_1, \cdots, n_m)$.   Then we have
\begin{equation}  \label{e1.1}
r(\X) \le \max \{ n_1, \cdots, n_m \}.
\end{equation}
\end{proposition}
{\bf Proof} Let $n=\max\{n_1,n_2,\cdots, n_m\}$ and $\A\in{\bf T}(n_1, \cdots, n_m)$ with a subtensor $\X$. Then  $r(\X)\leq r(\A)$ from Property 6. Together with  $r(\A)\leq n$ since $r$ is proper,  the result is arrived.

\qed

The matrix rank of an $n_1 \times n_2$ rectangular matrix is never greater than $\min \{ n_1, n_2 \}$.    From this proposition, we may think further to restrict the magnitude of the tensor rank.   For $m, n_1, \cdots, n_m \in {\bf N}$ with $m \ge 2$, we define submax$\{ n_1, \cdots, n_m \}$ as the second largest value of $n_1, \cdots, n_m$.

\begin{definition}
Suppose that $r$ is a tensor rank function.  We say that $r$ is a strongly proper tensor rank function if for any $m, n_1, \cdots, n_m \in {\bf N}$ with $m \ge 2$, and $\X \in {\bf T}(n_1, \cdots, n_m)$, we have
\begin{equation}  \label{e1.2}
r(\X) \le {\rm submax} \{ n_1, \cdots, n_m \}.
\end{equation}
\end{definition}

We cannot change submax$\{ n_1, \cdots, n_m \}$ in (\ref{e1.2}) to the third largest value of $n_1, \cdots, n_m$
as this violates Properties 1 and 3 of Definition \ref{d0}.

We will show that $r_*$ is strongly proper in the next section.

\section{CP Rank, Max-Tucker Rank and Submax-Tucker Rank}

As we stated in the introduction,  our motivation to introduce the axiom system for tensor ranks is to find some tensor ranks which have the max-full-rank-subtensor property.
The six properties of Definition \ref{d0} are not satisfied by some tensor ranks in the literature.
For example,
the tubal rank $r$ of third order tensors was introduced in \cite{KBHH13}.   For $\X \in T(n_1, n_2, n_3)$, $r(\X) \le \min \{ n_1, n_2 \}$.   Thus, it is not a tensor rank function even for third order tensors.     It is still very useful in applications \cite{YHHH16, ZSKA18, ZEAHK14, ZA17, ZLLZ18}.

However, the six properties of Definition \ref{d0} are satisfied by tensor ranks arising from two most important tensor decompositions -- the CP decomposition and the Tucker decomposition.

We now study the CP rank \cite{KB09}.

\begin{definition} \label{d1}
Suppose that $m, n_1, \cdots, n_m \in {\bf N}$ and $\X = (x_{i_1\cdots i_m}) \in {\bf T}(n_1, \cdots, n_m)$.   Suppose that there are $\aa^{(i, p)} \in \Re^{n_i}$ for $i=1, \cdots, m$ and $p = 1, \cdots, r$ such that
\begin{equation} \label{eq2.1}
\X = \sum_{p=1}^r \aa^{(1, p)} \circ \cdots \circ \aa^{(m, p)},
\end{equation}
then we say that $\X$ has a CP decomposition (\ref{eq2.1}).   The smallest integer $r$ such that (\ref{eq2.1}) holds is called the CP rank of $\X$, and denoted as $r_{CP}(\X)$.
\end{definition}

\begin{theorem}
The CP rank is a subadditive tensor rank function.   It is not a proper tensor rank function.
\end{theorem}
{\bf Proof}   We first show that the CP rank is a tensor rank function.   Properties 1, 3 and 4 hold clearly from the definition of the CP rank. Before we  show  Property 2, we can assert that $r_{CP}(\I_{m,n})\leq n$ for all $m,n \in {\bf N}$ with $m \ge 2$ since $\mathcal{I}_{m,n}=\sum_{i=1}^n \ee^i\circ \cdots \circ \ee^i$, where $\ee^i\in \Re^n$  with the unique nonzero entry  $e^i_i=1$. In the following, we show Property 2 by induction for $m$.   We fix $n$ here.

For $m=2$, $\mathcal{I}_{2,n}$ reduces to the $n \times n$ identity  matrix and hence Property 2 is true for such a case. Now we assume that $r_{CP}(\mathcal{I}_{m,n})=n$.
Then  we show
$r_{CP}(\mathcal{I}_{m+1,n})=n$.

Assume that $\mathcal{I}_{m+1,n}=\sum\limits_{p=1}^{r} \aa^{(1, p)} \circ \cdots \circ \aa^{(m+1, p)}$ with $r<n$. Then
\[\mathcal{I}_{m,n}=\mathcal{I}_{m+1,n}\cdot \ee \equiv \sum_{p=1}^r ((\ee)^T\aa^{(m+1,p)})\aa^{(1, p)} \circ \cdots \circ \aa^{(m, p)}.\]
Here, $\ee$ is the all one vector in $\Re^n$.
This indicates that $r_{CP}(\mathcal{I}_{m,n})<n$ since $r<n$. This contradicts the assumption that $r_{CP}(\mathcal{I}_{m,n})=n$ and  hence $r_{CP}(\mathcal{I}_{m+1,n})=n$.

  Hence, Property 2 holds.

For Property 5, we have that $\Y=\sum_{p=1}^r \aa^{(j_1, p)} \circ \dots \circ \aa^{(j_m, p)}$ if $\X=\sum_{p=1}^r \aa^{(1, p)} \circ \dots \circ \aa^{(m, p)}$ when $\Y$ is a permutation of $\X$ with $(j_1,\dots, j_m)=\sigma(1,2,\dots, m)$. Hence we have Property 5.

For property 6, assume that $\Y$ is a subtensor of $\X$. For $p=1,\dots, r$, let  $\X_p=\aa^{(1, p)} \circ \cdots \circ \aa^{(m, p)}$ and  $\Y_p$ be a subtensor of $\X_p$ by a similar way of $\Y$ from $\X$. Then we have that $\Y=\Y_1+\dots \Y_r$ and $r_{CP}(\Y)\leq r$ since $\Y_p$ are rank-one tensors for $p =1, \cdots, r$.
This means that $r_{CP}(\Y)\leq r_{CP}(\X)$.
%Furthermore, let $\Y$ be an essential subtensor of $\X$. Clearly, $r(\X)\geq r(\Y)$ since $\Y$ is a subtensor of $\X$.
%Now we assume that $r(\Y)<r(\X)$.  Suppose that $\Y=\sum_{p=1}^r \aa^{l_1,1}\circ \aa^{l_2,2}\circ\dots \circ \aa^{l_m,m}$ and let \[\bar \aa^{(i,p)}=\left\{\begin{array}{rll}&\aa^{l_i,i}, \quad &\mbox{if}~~(i,p)=(l_i,i)\in T_l,\\
%&\0,\quad &\mbox{otherwise,}\end{array}\right.\]
%where $T_l$ is the index set related to  subtensor $\Y$.
%Then $\X=\sum_{p=1}^r \bar \aa^{1,p}\circ \dots \bar \aa^{m,p}$ and hence $r(\X)\leq r$.
Hence %$r(\X)=r(\Y)$ and
Property 6 is satisfied.

Therefore, the CP rank is a tensor rank function.

Suppose that $\X,\Y\in {\bf T}(n_1, \dots, n_m)$ with $r_{CP}(\X)=r_1$ and $r_{CP}(\Y)=r_2$. Let
\[\X=\sum_{p=1}^{r_1} \aa^{(1, p)} \circ \cdots \circ \aa^{(m, p)}, \quad \Y=\sum_{q=1}^{r_2} \bb^{(1,q)}\circ \cdots \circ \bb^{(m,q)}.\]
It holds that
\[\X+\Y=\sum_{p=1}^{r_1} \aa^{(1, p)} \circ \cdots \circ \aa^{(m, p)}+\sum_{q=1}^{r_2} \bb^{(1,q)}\circ \cdots \circ \bb^{(m,q)}.\]

Hence, $r_{CP}(\X+\Y) \le r_1 + r_2 \equiv r_{CP}(\X) + r_{CP}(\Y)$.  This shows that it is subadditive.
By \cite{KB09}, the CP rank of a $9 \times 9 \times 9$ tensor given by Kruskal is between $18$ and $23$.
Thus, the CP rank is not a proper tensor rank function.
\qed

We now study the Tucker rank.  In some papers such as \cite{JYZ17}, the $n$-rank is called the Tucker rank.

Suppose that $m, n_1, \cdots, n_m \in {\bf N}$ and $\X = (x_{i_1\cdots i_m}) \in {\bf T}(n_1, \cdots, n_m)$.  We may unfold $\X$ to a matrix $X_{(j)} = (x_{i_j, i_1\cdots i_{j-1}i_{j+1} i_m}) \in \Re^{n_j \times n_1\cdots n_{j-1}n_{j+1}\cdots n_m}$ for $j = 1, \cdots, m$.   Denote $r_0(X_{(j)})$ as $r_j$ for $j = 1, \cdots, m$.  Then the vector $(r_1, \cdots, r_m)$ is called the $n$-rank of $\X$ \cite{KB09}.

The $n$-rank is a vector rank.   Hence it does not satisfy Definition \ref{d0}.  However, if we define
\begin{equation} \label{e3.1}
r = \max \{ r_1, \cdots, r_m \},
\end{equation}
then we have the following proposition.

\begin{theorem}\label{maxtuck}
The function $r$ defined by (\ref{e3.1}) is a proper, subadditive tensor rank function.   But it is not strongly proper.
\end{theorem}
{\bf Proof}   We first show that rank function $r$ defined by (\ref{e3.1}) is a tensor rank function. To see this, it suffices to show that Property 1-6 are all satisfied.

1) Suppose that $\X \in {\bf T}(n_1, \cdots, n_m)$ for $m, n_1, \cdots, n_m \in {\bf N}$ is a zero tensor.
Then $X_{(j)}$ are zero matrices for $j = 1, \cdots, m$.   This implies that $r_j(X_{(j)}) = 0$ for $j = 1, \cdots, m$.  By (\ref{e3.1}), we have $r(\X) = 0$.  On the other hand, assume that  $r(\X)=0$ for some $\X \in {\bf T}(n_1, \cdots, n_m)$ with $m, n_1, \cdots, n_m \in {\bf N}$.    This means that $r_i=0$ for $i=1, \cdots, m$, which means that  $X_{(i)}=0$ and hence $\X$ is a zero tensor.

Suppose that $r(\X)=1$, then $r_i(X_{(i)})
= 1$ for all $i=1,\cdots, m$. This can be seen as follows.  Assume that there exists $i_0$ such that $r_0(X_{(i_0)})=0$, then  $\X$ is a zero tensor since $X_{(i_0)}=0$.  From above analysis, $r(\X)=0$ if and only if $\X$ is a zero tensor. This contradicts with $r(\X)=1$.

Let $\X=\sum_{p=1}^{\bar r} \aa^{(1,p)}\circ \aa^{(2,p)}\dots \circ \aa^{(m,p)}$. Then
$X_{(1)}=\sum_{p=1}^{\bar r} \aa^{(1,p)}\circ (\aa^{(2,p)}\circ \dots \circ \aa^{(m,p)}).$ From $r_0(X_{(1)})=1$, we have that $\aa^{(1,p)} (p=1,\dots, \bar r)$ is rank-one.
From $X_{(2)}=\sum_{p=1}^r \aa^{(2,p)}\circ (\aa^{(1,p)}\circ \dots\circ \aa^{(m,p)})$ and $r_0(X_{(2)})=1$, we have that $\aa^{(2,p)}$ for all $p=1,\cdots, \bar r $ is also rank-one.

Similarly, we have that for any $i=1,\dots, m$,  $\aa^{(i,p)}$ $(p=1,2,\dots, \bar r)$ is rank-one.

Thus, $\X=\lambda \aa^{(1,1)}\circ \cdots\circ \aa^{(m,1)}$ for some $\lambda$ and hence $\X$
is a rank-one tensor.

Conversely, if $\X$ is a rank-one tensor, then $\X=\x^{(1)}\circ \cdots\circ \x^{(m)}$ for some nonzero vectors $\x^{(i)}\in \Re^{n_i}$. Then $X_{(i)}= \x^{(i)}\circ (\x^{(1)}\circ \dots\circ \x^{(m)})$ and $r_i(X_{(i)})=1$ for all $i=1,\ldots, m$. Thus $r(\X)=1$.

Based on the above analysis, Property 1 is satisfied.

2) Denote $\I \equiv \mathcal{I}_{m,n}.$ Then $I_{(i)}$ is a rectangular matrix which can be partitioned to an $n$-dimensional identity matrix and an $n \times (m-1)n$ zero matrix, for $i = 1, \cdots, m$, and hence $r_0(I_{(i)})=n$. Thus, $r(\I)=n$.

3) When $\X\in {\bf T}(n_1,n_2,1,\cdots, 1)$, we have $X_{(1)}\in \Re^{n_1\times n_2}$, $X_{(2)}=X_{(1)}^T\in \Re^{n_2\times n_1}$ and $X_{(i)}\in \Re^{1,n_1n_2}$ for any $i\geq 3$. Clearly, $r_1=r_2$ and hence $r(\X)=r_0(X_{(1)})=r_0(X_{(2)})$.

4) Suppose that $\X\in {\bf T}(n_1,n_2,\cdots, n_m)$. For any $\alpha \neq 0$, and any $i\in \{1,2,\cdots, m\}$, $(\alpha X)_{(i)}=\alpha X_{(i)}$ and hence $r_i(X_{(i)})=r_i((\alpha \X)_{(i)})$. Hence $r(\X)=r(\alpha\X)$.

5) Suppose that $\X\in {\bf T}(n_1,\cdots,n_m)$
 and $\Y$ is any permuted tensor of $\X$. Then $Y_{(i)}$ will be $X_{(j)}$ for some $j\in \{1,2,\cdots, m\}$. So $r_0(Y_{(i)})= r_0(X_{(j)})$.
Hence $r(Y)=\max\{r_0(Y_{(i)}): i=1,\cdots, m\}=\max\{r_0(X_{(j)}): j=1,\cdots, m\}=r(\X)$ and the result holds.

6) Suppose that $\Z$ is a subtensor of $\X$. Then for all $i=1,2,\cdots, m$, $Z_{(i)}$ will be a submatrix of $X_{(i)}$ and $r_0(Z_{(i)})\leq r_0(X_{(i)})$  since $r_0$
is the matrix rank. So $r(\Z)\leq r(\X)$. %Suppose that $\Y$ is an essential subtensor of $\X$, then $Y_{(i)}$ is an essential subtensor of $X_{(i)}$ and hence
%$r_i(Y_{(i)})=r_i(X_{(i)})$. So we can assert that $r(\Y)=r(\X)$.

Now we conclude that $r$ defined by (\ref{e3.1}) is a tensor rank function.

It is clear that such a tensor rank function $r$ is proper from its definition. Furthermore,  we have that such rank $r$ is also subadditive since  matrix rank is subadditive.

In addition, we consider $\X\in {\bf T}(3,2,2)$ with $X_{(1)}=[I;\ee]$ where $I$ is the identity matrix of three dimension. Hence $r(\X)=3>2={\rm submax}\{3,2,2\}$. Hence we conclude that such a tensor rank function is not strongly proper.

\qed

Thus, we call this tensor rank function the max-Tucker rank in this paper, and denote it as $r_{max}(\X)$ for any $\X \in {\bf T}$.

Note that the max-Tucker rank naturally arises from applications of the Tucker decomposition when people assume that $r_i \le r$ for $i = 1, \cdots, m$ and fix the value of $r$ \cite{CSZZC19, XY13}.  Then this means that tensors of max-Tucker ranks not greater than $r$ are used.   In the following, we introduce a new tensor rank function, which is also associated with the Tucker decomposition, but is different from the max-Tucker rank.  We may replace (\ref{e3.1}) by
\begin{equation} \label{e3.2}
r = {\rm submax} \{ r_1, \cdots, r_m \}.
\end{equation}
Then we have the following theorem.

\begin{theorem}
The function $r$ defined by (\ref{e3.2}) is a strongly proper tensor rank function.   But it is not subadditive.
\end{theorem}
{\bf Proof}
 We first show that function $r$ defined by (\ref{e3.2}) is a tensor rank function.
 It suffices to show that Property 1-6 are all satisfied.

1) Suppose that $\X \in {\bf T}(n_1, \cdots, n_m)$ for $m, n_1, \cdots, n_m \in {\bf N}$ is a zero tensor.
Then   $X_{(j)}$ are zero matrices for all $j = 1, \cdots, m$.   This implies that $r_0(X_{(j)}) = 0$, which means that $X_{(j)}=0$.   By (\ref{e3.2}), we have $r(\X) = 0$.  On the other hand, assume that  $r(\X)=0$ for some $\X \in {\bf T}(n_1, \cdots, n_m)$ with $m, n_1, \cdots, n_m \in {\bf N}$.    This means that for some $i\in\{1, \cdots, m\}$, $r_0(X_{(i)})=0$,  and hence $X_{(i)}=0$, $\X$ is a zero tensor.
Therefore, $\X$ is a zero tensor if and only if $r(\X)=0.$

Suppose that $r(\X)=1$. Then $\X$ is not a zero tensor and hence  then $r_0(X_{(i)})\geq 1$ for all $i=1,\cdots, m$. Since $r$ is defined by (\ref{e3.2}), there exists $i_1,i_2,\dots, i_{m-1}$ such that $r_{i_{j}}(\X)=r(X_{(i_j)})=1$. Without loss of generality, we assume that $i_j=j$ for $j=1,2,\dots, m-1$. Let $\X=\sum_{p=1}^{\bar r} \aa^{(1,p)}\circ \dots\circ \aa^{(m,p)}$.  Similar to discussion in proof of Theorem \ref{maxtuck},
 $\aa^{(j,p)}$ $(p=1,2,\dots, \bar r)$ is rank-one for all $j=1,\dots, m-1$. Thus
 \[\X= \aa^{(1,p)}\circ \dots\circ (\aa^{(m,1)}+\lambda_2\aa^{(m,2)}+\lambda_3 \aa^{(m,3)}+\cdots+\lambda_{\bar r}\aa^{(m,\bar r))},\]
for some $\lambda_2,\dots, \lambda_{\bar r}$. Clearly, such $\X$ is a rank-one tensor.

Conversely, if $\X$ is a rank-one tensor, then $\X=\x^{(1)}\circ \cdots\circ \x^{(m)}$ for some nonzero vectors $\x^{(i)}\in \Re^{n_i}$. Then $X_{(i)}= \x^{(i)}\circ (\x^{(1)}\circ \dots\circ \x^{(m)})$ and $r_i(X_{(i)})=1$ for all $i=1,\ldots, m$. Thus $r(\X)=1$.

Based on the above analysis, Property 1 is satisfied.

2) Denote $\I \equiv \mathcal{I}_{m,n}.$ Then $I_{(i)}$ is a rectangular matrix which can be partitioned to an $n$-dimensional identity matrix and an $n \times (m-1)n$ zero matrix, for $i = 1, \cdots, m$, and hence $r_0(I_{(i)})=n$. Thus, $r(\I)=n$.

3) When $\X\in {\bf T}(n_1,n_2,1,\cdots, 1)$, we have $X_{(1)}\in \Re^{n_1\times n_2}$, $X_{(2)}=X_{(1)}^T\in \Re^{n_2\times n_1}$ and $X_{(i)}\in \Re^{1\times n_1n_2}$ for any $i\geq 3$. Clearly, $r_0(X_{(1)})=r_0(X_{(2)})\geq 1$ and $r_0(X_{(i}))\leq 1$ when $i\geq 3$. Hence $r(\X)=r_1(X_{(1)})=r_2(X_{(2)})$ is the same as the matrix rank of the corresponding matrix.

4) Suppose that $\X\in {\bf T}(n_1,n_2,\cdots, n_m)$. For any $\alpha \neq 0$, and any $i\in \{1,2,\cdots, m\}$, $(\alpha \X)_{(i)}=\alpha \X_{(i)}$ and hence $r_0(X_{(i)})=r_0((\alpha \X)_{(i)})$. Hence $r(\X)=r(\alpha\X)$.

5) Suppose that $\X\in {\bf T}(n_1,\cdots,n_m)$
 and $\Y$ is any permuted tensor of $\X$. Then $Y_{(i)}$ will be $X_{(j)}$ for some $j\in \{1,2,\cdots, m\}$. So $r_0(Y_{(i)})= r_0(X_{(j)})$.
Hence $r(\Y)= {\rm submax} \{r_0(Y_{(i)}): i=1,\cdots, m\}= {\rm submax} \{r_0(X_{(j)}): j=1,\cdots, m\}=r(\X)$ and the result holds.

6) Suppose that $\Z$ is a subtensor of $\X$. Then for all $i=1,2,\cdots, m$, $Z_{(i)}$ will be a submatrix of $X_{(i)}$ and $r_0(Z_{(i)})\leq r_0(X_{(i)})$  since $r_i$
is matrix rank. So $r(\Z)\leq r(\X)$. %Suppose that $\Y$ is an essential subtensor of $\X$, then $Y_{(i)}$ is an essential subtensor of $X_{(i)}$ and hence
%$r_i(Y_{(i)})=r_i(X_{(i)})$. So we can assert that $r(\Y)=r(\X)$.

Now we conclude that $r$ defined by (\ref{e3.2}) is a tensor rank function.

 The strongly proper property of such a tensor rank function is clear and hence it suffices to show that it is not subadditive.

Let $\X = (x_{ijk}),\Y = (y_{ijk}),\Z = (z_{ijk})\in T(2n_1,2n_2,2n_3)$ with $\X=\Y+\Z$ and
\[y_{ijk}=0~~ \mbox{if}~~ i>n_1,~ j>n_2,~ k>n_3,\quad z_{pqs}=0 ~~\mbox{if}~~p\leq n_1,~q\leq n_2, s\leq n_3.\]
It is assumed that $n-rank(\Y)=(r_1,r_2,r_3)$, $n-rank(\Z)=(R_1,R_2,R_3)$ and $r_1>r_2>r_3$, $R_2>R_1>R_3$.
Then $X_{(i)}=Y_{(i)}+Z_{(i)}$ for $i=1,2,3$ and $r_0(X_{(i)})=r_0(Y_{(i)})+r_0(Z_{(i)})$. So
$r(\X)= {\rm submax}\{r_1+R_1,r_2+R_2, r_3+R_3\}> r_2+R_1$ since $r_1+R_1>r_2+R_1$ and $r_2+R_2>r_2+R_1$.
Therefore, we conclude that such a tensor rank function is not subadditive.

\qed

Thus, we call this tensor rank function the submax-Tucker rank in this paper, and denote it as $r_{sub}(\X)$ for any $\X \in {\bf T}$.

\begin{proposition}
We have  $r_{sub}(\Y)\le r_{max}(\Y)$ for any $\Y \in {\bf T}$ and $r_{sub}(\X) < r_{max}(\X)$ for some $\X \in {\bf T}$.  Thus,  $r_{max} \not = r_*$.    Furthermore, $r_*$ is strongly proper.
\end{proposition}
{\bf Proof}
Clearly, $r_{sub}(\Y)\le r_{max}(\Y)$ for any $\Y \in {\bf T}$. To see $r_{sub}(\X) < r_{max}(\X)$ for some $\X \in {\bf T}$, we consider the following counterexample $\X$.

Consider the tensor $\X\in {\bf T}(2,3,4)$ with its nonzeros entries $\X_{111}=\X_{122}=\X_{133}=\X_{214}=1$.
By observation, we have that $r_0(X_{(1)})=2, r_0(X_{(2)})=3$ and $r_0(X_{(3)})=4$, which implies that
$r_{sub}(\X)=3<4= r_{max}(\X)$ and the result is arrived here.

As $r_* \le r_{sub}$ and $r_{sub}$ is strongly proper, $r_*$ is also strongly proper.

\qed

We cannot replace submax $\{ r_1, \cdots, r_m \}$ in (\ref{e3.2}) by the third largest value in
$r_1, \cdots, r_m$, as this will violate Properties 1 and 3 of Definition \ref{d0}.

\section{Maximum Full Rank Subtensors}

 In this section, we introduce the concept of maximum full rank subtensors, and define the max-full-rank-subtensor property.

 We first define the full rank concept for a tensor rank function.    Recall that in matrix theory, there is the concept of full row (column) rank matrices.

\begin{definition} \label{d4.1}
Suppose that $r$ is a tensor rank function.  Let $m, n_1, \cdots, n_m \in {\bf N}$ with $m \ge 2$, and $\X \in {\bf T}(n_1, \cdots, n_m)$.   If we have
\begin{equation}  \label{e4.6}
r(\X) = n_p
\end{equation}
for some index $p$ satisfying $1 \le p \le m$,
then we say that $\X$ is of full $p$-row $r$ rank, or simply say that $\X$ is of full $r$ rank.
In particular,
zero tensors %and rank-one tensors
are regarded as of full $r$ rank.
\end{definition}

We then define the max-full-rank-subtensor property for a tensor rank function.

\begin{definition}  \label{d4.2}
Suppose that $r$ is a tensor rank function.
Let $\X \in {\bf T}$.   We call a subtensor $\Y$ of $\X$ a maximum full rank subtensor of $\X$ under $r$ if $\Y$ is of full $r$ rank, and $r(\Y)$ is the maximum for any such full rank subtensors of $\X$.
We say that $r$ is of the max-full-rank-subtensor property if for any $\X \in {\bf T}$,
$$r(\X) = r(\Y),$$
where $\Y$ is a maximum full rank subtensor of $\X$ under $r$.
\end{definition}

Now, the question is if there is a tensor rank function of the max-full-rank-subtensor property.   We have the following theorem.

\begin{theorem} \label{t4.3}
The max-Tucker rank function $r_{max}$ has the max-full-rank-subtensor property.

Furthermore, suppose that $m, n_1, \cdots, n_m \in {\bf N}$ and $\X = (x_{i_1\cdots i_m}) \in {\bf T}(n_1, \cdots, n_m)$, and $\Y$ is a maximum full rank subtensor of $\X$ under $r_{max}$.
Then there is an index $p$, $1 \le p \le m$, such that all the $q$th $p$-rows of $\X$, with $q$ in the mode $p$ index set $T_p$ of $\Y$, are linearly independent, $|T_p| = r_{max}(\Y) = r_{max}(\X)$, and any $p$-row of $\X$ is a linear combination of all the $q$th $p$-rows of $\X$ with $q \in T_p$.
\end{theorem}
{\bf Proof}   Let $m, n_1, \cdots, n_m \in {\bf N}$ with $m \ge 2$, and $\X \in {\bf T}(n_1, \cdots, n_m)$.  Assume that $r_0(X_{(i)}) = r_i$ for $i = 1, \cdots, m$. Without loss of generality, assume that $r_1 = \max \{ r_1, \cdots, r_m \}$.  By the properties of the matrix rank, we know that there is a set $T_1 = \{ k_1, \cdots, k_{r_1} \} \subset \{ 1, \cdots, n_1 \}$ such that $\Y = (y_{j_1\cdots j_m}) \in {\bf T}(r_1, n_2, \cdots, n_m)$ is a subtensor of $\X$, where
$$y_{j_1\cdots j_m} = x_{j_1\cdots j_m}$$
for $j_1 \in T_1$, $j_l = 1, \cdots, n_l$, $l = 2, \cdots, m$,
and $r_0(Y_{(1)})=r_1$.
Then $r_0(Y_{(l)}) \le r_0(X_{(l)}) \equiv r_l$ for $ l = 2, \cdots, m$.  This shows that $\Y$ is of full $r_{max}$ rank, and $r_{max}(\Y) = r_{max}(\X) = r_1$.    Hence, the max-Tucker rank function $r_{max}$ has the max-full-rank-subtensor property.

On the other hand, suppose that $m, n_1, \cdots, n_m \in {\bf N}$ and $\X = (x_{i_1\cdots i_m}) \in {\bf T}(n_1, \cdots, n_m)$, and $\Y$ is a maximum full rank subtensor of $\X$ under $r_{max}$.  By Definition \ref{d4.1}, there is an index $p$, $1 \le p \le m$, such that
$$r_{max}(\Y) = |T_p|,$$
where $T_p$ is the mode $p$ index set of $\Y$.   Denote $r_0(X_{(l)})$ and $r_0(Y_{(l)})$ as
$r_l(\X)$ and $r_l(\Y)$ respectively for $l = 1, \cdots, m$.   Then
$$r_l(\Y) \le r_l(\X)$$
for $l = 1, \cdots, m$.  We have
$$r_{max}(\Y) = |T_p| \le r_p(\Y) \le r_p(\X) \le r_{max}(\X) = r_{max}(\Y).$$
Hence,
$$r_{max}(\Y) = |T_p| = r_p(\Y) = r_p(\X) \le r_{max}(\X) = r_{max}(\Y).$$
This shows that all the $q$th $p$-rows of $\X$ with $q \in T_p$ are linearly independent by the definition of $r_p(\Y)$, and any $p$-row of $\X$ is a linear combination of all the $q$th $p$-rows of $\X$ with $q \in T_p$ by the definition of $r_p(\X)$.
\qed

The property of the maximum full rank subtensor $\Y$ of $\X$ under $r_{max}$, stated in Theorem \ref{t4.3}, extends the corresponding property of matrices to tensors.

\section{The Closure of a Proper Tensor Rank Function}

%With the discussion in the last two sections, we are now ready to study the sub-full-rank property in this section.  We first define the full rank concept for a strongly proper tensor rank function.

%\begin{definition}
%Suppose that $r$ is a strongly proper tensor rank function.  Let $m, n_1, \cdots, n_m \in {\bf N}$ with %$m \ge 2$, and $\X \in {\bf T}(n_1, \cdots, n_m)$.   If we have
%\begin{equation}  \label{e1.3}
%r(\X) = {\rm submax} \{ n_1, \cdots, n_m \},
%\end{equation}
%then we say that $\X$ is of full $r$ rank.
%\end{definition}

%We then define the sub-full-rank property for a strongly proper tensor rank function.

%\begin{definition}  \label{d4.2}
%Suppose that $r$ is a strongly proper tensor rank function.  We say that $r$ is of the sub-full-rank property if for any $m, n_1, \cdots, n_m \in {\bf N}$ with $m \ge 2$, and $\X \in {\bf T}(n_1, \cdots, n_m)$, either $\X$ is of full $r$ rank, or $\X$ has a subtensor $\Y$ such that $r(\X) = r(\Y)$ and $\Y$ is of full $r$ rank.   In particular,
%zero tensors and rank-one tensors are regarded as full of $r$ rank.
%\end{definition}

Theorem \ref{t4.3} says that the max-Tucker rank function $r_{max}$ has the max-full-rank-subtensor property.  Is there any other tensor rank function which also has this property?    The full rank concept is only suitable for proper tensor rank functions.   Thus, the CP rank $r_{CP}$ is out of question.   Does the submax-Tucker rank function $r_{sub}$ have the max-full-rank-subtensor property?   At this moment, we do not know the answer to this question.   However, we show that for any proper tensor rank function, we may also derive another proper tensor rank function, which has
the max-full-rank-subtensor property. To do this, we introduce the concept of the closure of a proper tensor rank function.

\begin{definition}  \label{d4.4}
Suppose that $r : {\bf T} \to {\bf Z_+}$ is a proper tensor rank function.
We may define $\bar r : {\bf T} \to {\bf Z_+}$ as the closure of $r$ by
$$\bar r(\X) = \max \{ r(\Y) : \Y \ {\rm is \ a \ subtensor \ of}\ \X, \ {\rm and \ of \ full}\ r \ {\rm rank} \},$$
for any $\X \in {\bf T}$.
\end{definition}

We have the following theorems.

\begin{theorem} \label{t5.2}
The closure $\bar r$ of a proper tensor rank function $r$ is also a
proper tensor rank function.  We have $\bar r \le r$.
A proper tensor rank function $r$ has the max-full-rank-subtensor property if and only if $\bar r = r$.
\end{theorem}
{\bf Proof}   Let $\X$ be a zero tensor.   By Definition \ref{d4.1}, $\X$ is of full $r$ rank.   By
Definition \ref{d4.4}, $\bar r(\X) = 0$.    Let $\X$ be a nonzero rank-one tensor.   Then $\X$ has a one-entry nonzero subtensor $\Y$ and $r(\Y) = 1$.   By Definition \ref{d4.1}, $\Y$ is of full $r$ rank.
Thus, $r(\Y) = 1$.   This shows $\bar r(\X) \ge 1$.    By Property 6 of Definition \ref{d0}, for any subtensor $\Z$ of $\X$, $r(\Z) \le r(\X) \le 1$.  This shows that $\bar r(\X) \le 1$.   Hence, $\bar r(\X) = 1$.  This shows that $\bar r$ satisfies Property 1 of Definition \ref{d0}.

By Property 2 of Definition \ref{d0}, for $m, n \in {\bf N}$ with $m \ge 2$, $r(\I_{m, n}) = n$.   By Definition \ref{d4.1}, $I_{m, n}$ is of full $r$ rank.   By Definition \ref{d4.4},  $\bar r(\I_{m, n}) = n$.   Thus, $\bar r$ satisfies Property 2 of Definition \ref{d0}.

Assume that $n_1, n_2 \in {\bf N}$, $\X \in {\bf T}(n_1, n_2, 1, 1, \cdots, 1)$.    Let $M$ be the corresponding $n_1 \times n_2$ matrix.   Denote that $r_0(M)=r_M$.   Then there is an
$r_M \times r_M$ submatrix $\bar M$ of $M$ such that the matrix rank of $\bar M$ is $r_M = r(\X)$.  Furthermore, there is a subtensor $\Y$ of $\X$ such that $\Y \in {\bf T}(r_M, r_M, 1, 1, \cdots, 1)$, and the corresponding $r_M \times r_M$ matrix is $\bar M$.   We now see that $\Y$ is of full $r$ rank.  This
shows that $\bar r(\X) \ge r(\Y) = r_M = r(\X)$.   Thus, $\bar r(\X) = r_M$.  This shows that $\bar r$ satisfies Property 3 of Definition \ref{d0}.

For any $\alpha\neq 0$, $r(\alpha\X)=r(\X)$ and hence
\[\begin{array}{rl}\bar r(\alpha \X)&=\max\{r(\alpha \Y): \Y \ {\rm is \ a \ subtensor \ of}\ \X, \ {\rm and \ of \ full}\ r \ {\rm rank} \}\\
&=
\max\{r( \Y): \Y \ {\rm is \ a \ subtensor \ of}\ \X, \ {\rm and \ of \ full}\ r \ {\rm rank} \}=\bar r(\X).\end{array}\]
This means that Property 4 is satisfied. Similarly, we have Property 5 for $\bar r$.

For any subtensor $\Z$ of $\X$, we have
\[\begin{array}{rl}\bar r(\Z)=&\max\{r(\Y): \Y \ {\rm is \ a \ subtensor \ of}\ \Z, \ {\rm and \ of \ full}\ r \ {\rm rank} \}\\
\leq &\max\{r(\Y): \Y \ {\rm is \ a \ subtensor \ of}\ \X, \ {\rm and \ of \ full}\ r \ {\rm rank} \}=\bar r(\X).
\end{array}
\]
%Furthermore, suppose that $\Z$ is any essential subtensor of $\X$, then $r(\Z)=r(\X)$ and all entries in $\X$  but not in  $\Z$ are zeros. %Then
%\[\bar r(\Z)=\max\{r(\Y): \Y \ {\rm is \ a \ subtensor \ of}\ \Z, \ {\rm and \ of \ full}\ r \ {\rm rank} \}.\]
%Let \[\bar r(\X)=\max\{r(\Y): \Y \ {\rm is \ a\ subtensor \ of}\ \X, \ {\rm and \ of \ full}\ r \ {\rm rank} \}
%:= r( \Y),
%\] where  $\Y$ is a subtensor of $\X$ of full $r$ rank.
%If $\Y=\Z$, then $\bar r(\Z)=r(\Y)=\bar r(\X)$.
%Otherwise,
%let $\W$ be the intersection subtensor of $\Z$ and $\Y$.   Note that the entries that are in $\Y$ but not in $\W$ must be zeros since such entries must not be in $\Y$. This means that $\W$ is an essential subtensor of $\Y$ and hence $r(\W)=r(\Y)$. Since $r$ is proper, $\W$ is a full $r$ rank subtensor of $\Z$, $\Y$ and $\X$, respectively. Then
 %\[  \bar r(\Z)\geq r(\W)=r(\Y)=\bar r(\X)\geq \bar r(\Z).\]
 %This means that $\bar r(\Z)=\bar r(\X)$ and hence
Hence, Property 6 of Definition \ref{d0} is satisfied by $\bar r$.

Clearly, $\bar r$ is proper since $r$ is proper. So we can assert that
 $\bar r$ is a proper tensor rank function.

By Property 6 of $r$ and the definition of $\bar r$, we have $\bar r \le r$.

 Now we show the last assertion that $r$ has the max-full-rank-subtensor property if and only if $\bar r=r$.

``$\Rightarrow$'' It suffices to show that $\bar r\geq r$. If $\X$ is of full $r$ rank, then $\bar r(\X)=r(\X)$.  Otherwise, there exists a full $r$ rank subtensor $\Y$ of $\X$ such that $r(\Y)=r(\X)$. Since $\bar r(\X)\geq r(\Y)$ from definition, $\bar r (\X)\geq r(\X)$. Together with $\bar r (\X)\leq r(\X)$, we have $\bar r(\X)=r(\X).$

``$\Leftarrow$'' From $r(\X)=\bar r(\X)$ for any $\X$, we have that there exists a full $r$ rank subtensor $\Y$ of $\X$ such that
$\bar r(\X)=r(\Y)$. So $r(\Y)=r(\X)$ and hence $r$ has the max-full-rank-subtensor property from the arbitrariness of $\X$.

 The conclusion holds.

\qed

\begin{theorem}
Suppose that $r$ is a proper tensor rank function, and $\bar r$ is its closure.   Then $\bar r$ has the max-full-rank-subtensor property.
\end{theorem}
{\bf Proof}   Let $\X \in {\bf T}$ and $\Y$ be a maximum full rank subtensor of $\X$ under $r$.
By Definition \ref{d4.4}, $\bar r(\X) = r(\Y)$, and $\bar r(\Y) = r(\Y)$.  Then $\Y$ is also a maximum full rank subtensor of $\X$ under $\bar r$,  and we have
$$\bar r(\Y) = \bar r(\X).$$
This shows that $\bar r$ has the max-full-rank-subtensor property.
\qed

Then we are now able to show that the smallest tensor rank function $r_*$ is such a strongly proper tensor rank function.

\begin{corollary}
The smallest tensor rank function $r_*$ has the max-full-rank-subtensor property.
\end{corollary}
{\bf Proof}   Let $r_{**}$ be the closure of $r_*$.  Since, $r_* \le r_{**} \le r_*$, we have $r_* = r_{**}$.   Hence, $r_*$ has the max-full-rank-subtensor property.
\qed

Is $r_*$ equal to the submax-Tucker rank function $r_{sub}$ or its closure $\bar r_{sub}$?   This leaves as a further research question.

\section{An Application of The Submax-Tucker Rank}

In Section 3, we introduced a new tensor rank function, the submax-Tucker rank function, which is associated with the Tucker decomposition, but is different from the max-Tucker rank.
According to our theoretical analysis, the submax-Tucker rank is strongly proper.   Comparing with the CP rank and the max-Tucker rank, it is smaller in general.   Thus, the submax-Tucker rank may be a good choice for low rank tensor approximation and tensor completion.   We now present an application of the submax-Tucker rank.

Suppose that we have a data tensor $\M \in {\bf T}(n_1, n_2, \cdots, n_m)$.   Assume that $n_1 >> n_i$ for $i=2, \cdots, m$.   Then we may approximate $\M$ by $\X \in {\bf T}(\bar n, r, \cdots, r)$, where
$n_1 \ge \bar n \ge n_i$ for $i = 2, \cdots, m$, and $r \le \max \{ n_2, \cdots, n_m \}$. ,   For example, in \cite{ZZXC15}, for the internet traffic data tensor Abilene $\M$ \cite{Abilene}, we have $n_1 = 1008$, which is the number of time intervals, $n_2 = n_3 = 11$ is the number of the origin-destination nodes of the internet traffic dataset.  We may use the Tucker decomposition \cite{KB09}
$$\X = \D \times_1 A \times_2 B \times_3 C$$
to approximate $\M$.  Here, $\D$ is the Tucker core tensor of dimension $r_1 \times r_2 \times r_3$.   Factor matrices  $A, B$ and $C$ are of dimensions $n_1 \times r_1$, $n_2 \times r_2$ and $n_3 \times r_3$, respectively.   The operations $\times_i$ are mode $i$ product \cite{KB09}.   A usual practice is to fix $r$ and assume that $r_i \le r$ for $i = 1, 2, 3$ \cite{CSZZC19, XY13}.   This means to approximate the data tensor $\M$ by a tensor $\X$ of the max-Tucker rank not greater than $r$.   Then the range $r$ is $1 \le r \le 1008$, though we always have $r_2 \le 11$ and $r_3 \le 11$.   The range of $r$ is quite large.   If we use a tensor $\X$ of the submax-Tucker rank not greater than $r$ to approximate $\M$, then the range of $r$ is $1 \le r \le 11$, we may let, say, $r_1 \le \bar n = 30$, $r_2 \le r$ and $r_3 \le r$, by fixing $r$.   This provides a good choice of the range of $\X$ to approximate $\M$.

%In mathematics, there are examples, in which not the most extreme value, but the second most extreme value is the most useful.    For example, in spectral graph theory \cite{BH11, C97}, the second smallest eigenvalue of the Lapalacian tensor is called algebraic connectivity, and is the most useful.

\begin{figure}
  \centering
  \includegraphics[width=0.88\textwidth]{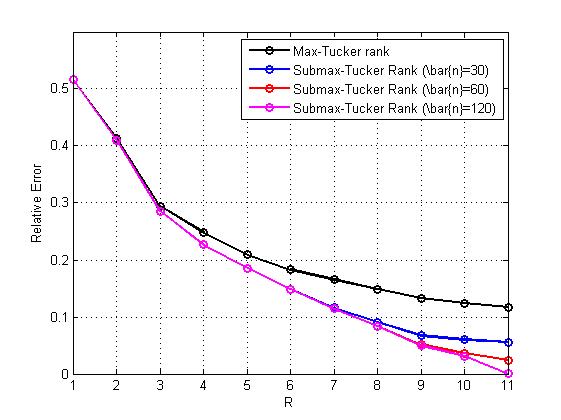}\\
  \caption{Comparison between the max-Tucker rank and the submax-Tucker Ranks.}\label{Abilene}
\end{figure}

For example, we consider the internet traffic tensor $\X\in{\bf T}(1008,11,11)$. We compare four kinds of Tucker decompositions of $\X$. (I) Tucker decomposition with the max-Tucker rank $r$, i.e., the core tensor $\D_1\in{\bf T}(r,r,r)$. (II--IV) Tucker decomposition with the submax-Tucker rank $r$ and $\bar{n}=30,60,120$, i.e., the core tensor $\D_1\in{\bf T}(\bar{n},r,r)$, respectively. For each decomposition $\widetilde{\X}$, we calculate the relative error
$$ \mathrm{Relative~error} := \frac{\|\widetilde{\X}-\X\|_F}{\|\X\|_F}.  $$
Using the Tensor Toolbox, we illustrate results in Figure \ref{Abilene} for $r$ ranging from $1$ to $11$. Obviously, we see that relative~errors corresponding to submax-Tucker rank is smaller than the relative~error of the max-Tucker rank case.

\section{Final Remarks}

In this paper, we extended the maximum full rank subtensor concept and  the max-full-rank-submatrix property to tensors.   We proved that the max-Tucker rank function, the smallest tensor rank function, and the closure of any proper tensor rank function have the max-full-rank-subtensor property.    These show that the maximum full rank subtensor concept and  the max-full-rank-subtensor property should be an important part for the tensor rank theory.   Some questions remain for further research.

The axiom system for tensor ranks is also an exploration.
The six properties of Definition \ref{d0} may be further modified.   But it may be a worthwhile research direction to study tensor ranks with some appropriate axiom systems.

%(3) We studied tensor rank functions associated with the CP decomposition and the Tucker decomposition.   In particular, we introduced the submax-Tucker rank.    It is possible that $r_*$ is the submax-Tucker rank function or its closure.   This will be a further research topic.  In this paper, we only
%discussed an application of the submax-Tucker rank.

%There are some other tensor ranks in the literature, which do not satisfy our axiom systems, but are still very useful.

For low rank tensor approximation, the concept of border rank \cite{KB09} is useful.   Can the concept of border rank also be accommodated by the tensor rank axiom system?   This may also be an interesting further research topic.

\bigskip

%{\bf Acknowledgment}  We are thankful to Haibin Chen and Ziyan Luo for their comments.

\end{document}